\documentclass[preprint,3p]{elsarticle}
\usepackage{amsmath}
\usepackage{amssymb}
\newtheorem{thm}{Theorem}[section]
\newtheorem{lem}[thm]{Lemma}
\newtheorem{cor}[thm]{Corollary}
\newtheorem{prop}{Proposition}[section]
\newtheorem{rem}{Remark}[section]

\newdefinition{rmk}{Remark}[section]
\newproof{pf}{Proof}
\newproof{pot}{Proof of Theorem \ref{thm2}}
\newproof{poot}{Proof of Corollary \ref{co1}}
\numberwithin{equation}{section}
\newdefinition{ex}{Example}[section]
\journal{}
\begin{document}
\begin{frontmatter}

\title{  Regularity of  extremal solutions of semilinaer fourth-order elliptic problems with general nonlinearities  }
\author{A. AGHAJANI}
\ead{aghajani@iust.ac.ir}.

\address{School of Mathematics, Iran University of Science and Technology, Narmak, Tehran, Iran.}
\address{School of Mathematics, Institute for Research in Fundamental Sciences (IPM), P.O.Box: 19395-5746, Tehran, Iran.}

\begin{abstract}
We consider the fourth order problem $\Delta^{2}u=\lambda f(u)$ on a general bounded domain $\Omega$ in $R^{n}$ with the Navier boundary condition $u=\Delta u=0$ on $\partial \Omega$. Here, $\lambda$ is a positive parameter and $ f:[0,a_{f}) \rightarrow \Bbb{R}_{+} $ $ (0 < a_{f} \leqslant \infty)$ is a smooth, increasing, convex nonlinearity such that $ f(0) > 0 $ and which blows up at $ a_{f} $. Let
$$0<\tau_{-}:=\liminf_{t\rightarrow a_{f}} \frac{f(t)f''(t)}{f'(t)^{2}}\leq \tau_{+}:=\limsup_{t\rightarrow a_{f}} \frac{f(t)f''(t)}{f'(t)^{2}}<2.$$
We show that if $u_{m}$ is a sequence of semistable solutions   correspond to $\lambda_{m}$ satisfy the stability inequality
 $$ \sqrt{\lambda_{m}}\int_{\Omega}\sqrt{f'(u_{m})}\phi^{2}dx\leq \int_{\Omega}|\nabla\phi|^{2}dx, ~~\text{for all}~\phi\in H^{1}_{0}(\Omega),$$
 then $\sup_{m} ||u_{m}||_{L^{\infty}(\Omega)}<a_{f}$ for $n< \frac{4\alpha_{*}(2-\tau_{+})+2\tau_{+}}{\tau_{+}}\max \{1, \tau_{+}\},$ where $\alpha^{*}$ is the largest root of the equation
$$(2-\tau_{-})^{2} \alpha^{4}- 8(2-\tau_{+})\alpha^{2}+4(4-3\tau_{+})\alpha-4(1-\tau_{+})=0.$$
In particular, if $\tau_{-}=\tau_{+}:=\tau$, then     $\sup_{m} ||u_{m}||_{L^{\infty}(\Omega)}<a_{f}$ for $n\leq12$ when $\tau\leq 1$,  and  for $n\leq7$ when $\tau\leq 1.57863$. These estimates lead to the regularity of the corresponding extremal solution $u^{*}(x)=\lim_{\lambda\uparrow\lambda^{*}}u_{\lambda}(x),$  where $\lambda^*$ is the extremal parameter of the eigenvalue problem.
\end{abstract}

\end{frontmatter}
\section{Introduction and main results}
In this article, we consider the  problem

\begin{equation*}
\left\{\begin{array}{ll} \Delta^{2} u= \lambda f(u)& {\rm }\ x\in \Omega,\\~~u=\Delta u=0& {\rm }\ x\in \partial \Omega,~~~~~~~~~~(N_{\lambda})
\end{array}\right.
\end{equation*}
where $\Omega\subset R^{n}$ is a smooth bounded domain, $n\geq1$, $\lambda >0$ is a real parameter, and the nonlinearity $f$ satisfies\\
$(H)$ $ f:[0,a_{f}) \rightarrow \Bbb{R}_{+} $ $ (0 < a_{f} \leqslant \infty)$ is a smooth, increasing, convex function such that $ f(0) > 0 $ and $\lim_{t\rightarrow a_{f}}f(t)=\infty.$ Also, when $a_{f}=\infty$ we assume that $f$ is superlinear, i.e., $\lim_{t\rightarrow\infty}\frac{f(t)}{t}=\infty$. \\
We call the nonlinearity $f$ \emph{regular}  if $a_{f}=\infty$ and \emph{singular} when $a_{f}<\infty$.

By a semistable solution  of  $N_{\lambda}$ we mean a solution $u$ satisfies
\begin{equation}
\int_{\Omega}(\Delta \varphi)^{2}-\int_{\Omega}\lambda f(u)'\varphi^{2}\geq0,~~~\varphi\in H^{2}(\Omega)\cap H^{1}_{0}(\Omega).
\end{equation}
Also, we say that a smooth solution $u$ of $N_{\lambda}$  is minimal provided $u\leq v$ $a.e.$ in $\Omega$ for any solution $v$ of $N_{\lambda}$ (see \cite{CGE,CG}.\\

When $f$ satisfies ($H$) is a regular, or $f(t)=(1-t)^{-p}$ ($p>1$), it is well known \cite{BG,CaG,GW} that there exists a finite positive extremal parameter $\lambda^{*}>0$ depending on $f$ and $\Omega$ such that for any $0<\lambda<\lambda^{*}$, problem ($N_{\lambda}$) has a minimal smooth solution $u_{\lambda}$, which is semistable and unique among  the semistable solutions,  while no solution exists for $\lambda\geq\lambda^{*}$. The function $\lambda\rightarrow u_{\lambda}$ is strictly increasing on $(0,\lambda_{*})$, the increasing pointwise limit $u^{*}(x)=\lim_{\lambda\uparrow\lambda^{*}}u_{\lambda}(x)$ is called the extremal solution. For $0<\lambda<\lambda_{*}$ the minimal solution $u_{\lambda}$   of problem $(N_{\lambda})$ satisfies the following stability inequality, for the proof see Corollary 1 in \cite{CG} or Lemma 6.1 in \cite{Dup1},
 \begin{equation}
 \sqrt{\lambda}\int_{\Omega}\sqrt{f'(u_{\lambda})}\phi^{2}dx\leq \int_{\Omega}|\nabla\phi|^{2}dx,
\end{equation}
for all $\phi\in H^{1}_{0}(\Omega)$.\\


The regularity and properties of the extremal solutions have been studied extensively in the literature [2-12,15,19] and it is shown that it depends strongly on the dimension $n$, domain $\Omega$ and nonlinearity $f$. \\
 Cowan,  Esposito and  Ghoussoub in \cite{CGE} showed that for general nonlinearity $f$ satisfies (H), $u^{*}$ is bounded for $n\leq 5$. When  $f(u)=e^{u}$, in \cite{CGE} it is shown that $u^{*}$ is bounded for $n\leq 8$.  This result improved by Cowan and Ghoussoub to $n\leq 10$ in \cite{CG}, and by Dupaigne,  Ghergu and  Warnault in \cite{Dup1} to $n\leq 12$ which is the optimal dimension as we know on the unit ball $u^{*}$ is bounded if and only if $n\leq 12$. As we shall see, in this paper we prove the same for a large class of nonlinearities including $e^{u}$. When $f(u)=(1+u)^{p}$ ($p>1$) in \cite{CGE} it is proved that $u^{*}$ is bounded if $n<\frac{8p}{p-1}$ that improved in \cite{CG} for  to $n<4h(p)>\frac{8p}{p-1} $ (for the definition of $h(p)$ which is a decreasing function on $(1,\infty)$ see \cite{CG}) with $\lim_{p\rightarrow\infty}4h(p)\approx10.718$. Recently, Hajlaoui, Harrabi and  Ye in \cite{HHY}  improved this result by showing that   $u^{*}$ is bounded  for any $p>1$ and $n\leq 12$.\\

 For the singular nonlinearity $f(u)=(1-u)^{-p}$ ($p>1$), in \cite{CGE} it is proved that $\sup_{\Omega}u^{*}<1$  if $n\leq\frac{8p}{p+1}$. In particular, when $p=2$, $u^{*}$ is bounded away from $1$ for $n\leq5$. The later result (and also the general case $1<p\neq 3$) is improved in \cite{CG}  to $n\leq6$, and further improved by Guo and Wei in \cite{GuW} to $n\leq 7$. However, for $p=2$ the expected optimal dimension is $n=8$, holds on the ball, see \cite{Amir}.\\


By imposing extra assumptions on the general nonlinearity $f$ satisfies (H), the authors in  \cite{CGE} obtained more regularity results in higher dimensions on general domains. Let $f$  satisfy ($H$) and define
\begin{equation}
\tau_{-}:=\liminf_{t\rightarrow a_{f}} \frac{f(t)f''(t)}{f'(t)^{2}}\leq \tau_{+}:=\limsup_{t\rightarrow a_{f}} \frac{f(t)f''(t)}{f'(t)^{2}}.
\end{equation}\\

 In \cite{CGE} the authors also show that for a regular and superlinear nonlinearity $f$ with $\tau_{-}>0$, $u^{*}$ is bounded for $n\leq 7$ (see \cite{CGE}, Theorem 4.1). As we shall see here in Corollary 2.4, with a minor change in their proof, the same holds with a  weaker condition. Also, they showed that if $\tau_{+}<\infty$ then $u^{*}$ is bounded for $n< \frac{8}{\tau_{+}}$, see Theorem 5.1 in \cite{CGE}.\\


The main results of this paper are as follows.\\
\begin{thm}
 Let $f$  satisfy $(H)$ with $0<\tau_{-}\leq\tau_{+}<2$, and $\Omega$ an arbitrary bounded smooth domain. Also, let  $u_{m}$ be a sequence of semistable solutions of $(N_{\lambda_{m}})$ satisfy the stability inequality (1.2). Then  $\sup_{m} ||u_{m}||_{L^{\infty}(\Omega)}<a_{f}$ for
\begin{equation}
n<N(f):=\left\{\begin{array}{ll} \frac{4\alpha_{*}(2-\tau_{+})+2\tau_{+}}{\tau_{+}}\max\{1,\tau_{+}\}& {\rm }\ f~is~regular,\\~~\frac{4\alpha_{*}(2-\tau_{+})+2\tau_{+}}{\tau_{+}}& {\rm }\ f~is~singular,
\end{array}\right.
\end{equation}
where $\alpha_{*}>1$ denotes  the largest root of the polynomial
\begin{equation}
P_{f}(\alpha,\tau_{-},\tau_{+}):=(2-\tau_{-})^{2} \alpha^{4}- 8(2-\tau_{+})\alpha^{2}+4(4-3\tau_{+})\alpha-4(1-\tau_{+}).
\end{equation}
As a consequence,\\
if $\tau_{-}=\tau_{+}:=\tau$, then     $\sup_{m} ||u_{m}||_{L^{\infty}(\Omega)}<a_{f}$ for $n\leq12$ when $\tau\leq 1$,  and  for $n\leq7$ when $\tau\leq 1.57863$.
\end{thm}

\begin{cor}
Let $f$  satisfy $(H)$ be a regular nonlinearity with $0<\tau_{-}\leq\tau_{+}<2$ and $\Omega$ an arbitrary bounded smooth domain. Let $u^{*}$ be the extremal solution of problem $(N_{\lambda})$. Then $u^{*}\in L^{\infty}(\Omega)$ for
$$n<\frac{4\alpha_{*}(2-\tau_{+})+2\tau_{+}}{\tau_{+}}\max \{1, \tau_{+}\}.$$
In particular, if $\tau_{-}=\tau_{+}$ then $u^{*}\in L^{\infty}(\Omega)$ for $n\leq12$.
\end{cor}

For example consider problem (1.1) with  $f(u)=e^{u}$ or $e^{u^{\alpha}}$ ($\alpha>0$), then $\tau_{+}=\tau_{-}=1$, hence by Theorem 1.1, $u^{*}\in L^{\infty}(\Omega)$ for $n\leq 12$. The same is true for $f(u)=(1+u)^{p}$ ($p>1$) as in this case  we have  $\tau_{+}=\tau_{-}=\frac{p-1}{p}$. More precisely we have $u^{*}\in L^{\infty}(\Omega)$ for $n<\frac{4(p+1)}{p-1}\alpha_{*}+2$
where $\alpha_{*}$ denotes  the largest root of the polynomial
\begin{equation}
P_{f}(\alpha):= (p+1)^{2} \alpha^{4}- 8p(p+1)\alpha^{2}+4p(p+3)\alpha-4p.
\end{equation}
This is exactly the same as the result obtained by Hajlaoui-Harrabi-Ye in \cite{HHY}.

Now consider problem (1.1) with  the singular nonlinearity $f(u)=(1-u)^{-p}$ $(p>1)$  and $\Omega$ an arbitrary bounded smooth domain. Then from the fact that $\tau_{+}=\tau_{-}=\frac{p+1}{p}$ and Theorem 1.1, we get $||u^{*}||_{L^{\infty}(\Omega)}<1$ for
$n<4\alpha_{*}+2,$ where $\alpha_{*}$ denotes  the largest root of the polynomial
\begin{equation}
P_{f}(\alpha):=  \alpha^{4}- 8\frac{p(p-1)}{(p+1)^{2}}\alpha^{2}+4\frac{p(p-1)(p-3)}{(p+1)^{3}}\alpha+4\frac{p(p-1)^{2}}{(p+1)^{4}}.
\end{equation}
This results coincides with that of Guo-Wei \cite{GW}. In particular, when $p>1.72822$ then $||u^{*}||_{L^{\infty}(\Omega)}<1$ for $n\leq 7$. Also, when $p>2.2609$  the same is true for $n\leq 8$.

\section{Preliminaries and Auxiliary Results}
The following standard regularity result is taken from \cite{CSS}, for the proof see Theorem 3 of \cite{Ser}.
\begin{prop}
Let $u\in H^{1}_{0}(\Omega)$ be a weak solution of
\begin{equation}
\left\{\begin{array}{ll} \Delta u+c(x)u=g(x)& {\rm }\ x\in \Omega,\\~~u=0& {\rm }\ x\in \partial \Omega,
\end{array}\right.
\end{equation}
with $c,g\in L^{q}(\Omega)$ for some $q> \frac{n}{2}$.
Then there exists a positive constant $ C$  independent of $u$ such that:
\begin{equation}
||u||_{L^{\infty}(\Omega)}\leq C(|u||_{L^{1}(\Omega)}+|g||_{L^{p}(\Omega)}).
\end{equation}
\end{prop}

Consider problem ($N_{\lambda}$). By the elliptic regularity we know that, if for some $q\geq 1$ we have $||f(u_{\lambda})||_{L^{q}(\Omega)}\leq C$, where $C$ is a constant independent of $\lambda$, then $u^{*}$ is bounded, (hence smooth when $f$ is regular), whenever $n<4q$. Using the above proposition we show that,  a similar result holds (for regular or singular nonlinearity) if $f'(u_{\lambda})$ is uniformly bounded in $L^{q}(\Omega)$. For the proof we  need the following two lemmas, the first one gives pointwise estimate
on $\Delta u$ for a solution $u$ of problem ($N_{\lambda}$), for the proof see \cite{CGE}.\\

Define the functions $F,g, \tilde{f}:[0,a_{f})\rightarrow R$ as
\begin{equation}
F(t)=\int_{0}^{t}f(s)ds,~~g(t)=\sqrt{2}(F(t)-t)^{\frac{1}{2}}~~\text{and}~~ \tilde{f}(t)=f(t)-f(0),~~0\leq t< a_{f}.
\end{equation}

\begin{lem}(Lemma 3.2 \cite{CGE})
Let $u$ be a solution of problem (1.1). Then
$$-\Delta u\geq \sqrt{\lambda} g(u),~~\text{in}~\Omega.$$
\end{lem}
 \begin{lem}
 Let $u$ be a semistable solution of problem problem $(N_{\lambda})$ with $f$ satisfy (H). Then
 $$\int_{\Omega}-\Delta u dx<C,$$
 where $C$ is a constant independent of $u$.
\end{lem}
\begin{pf}
Let $\psi$ be the unique positive smooth function such that
\begin{equation*}
\left\{\begin{array}{ll} -\Delta \psi= 1& {\rm }\ x\in \Omega,\\~~\psi=0& {\rm }\ x\in \partial \Omega,
\end{array}\right.
\end{equation*}
Let $u$ be a semistable solution of problem ($N_{\lambda}$). By multiplying the  equation $\Delta^{2}u=\lambda f(u)$ in $\psi$ and then an integration we get (using Green's formula)
$$\lambda\int_{\Omega}\psi(x)f(u)dx=\int_{\Omega}\psi(x)\Delta^{2}u~dx=\int_{\Omega}\Delta\psi(x)\Delta u~dx=\int_{\Omega}-\Delta u~dx.$$
This gives that
$$\int_{\Omega}-\Delta u~dx\leq \lambda \max_{\Omega}\psi(x)\int_{\Omega}f(u)dx.$$
The inequality above and the uniform $L^{1}(\Omega)$ boundedness of $f(u)$ for semistable solutions (proved in Lemma 3.5 in \cite{CGE}) gives the desired result. $\blacksquare$
\end{pf}
In the sequel we will frequently use the following simple lemma.
\begin{lem}
Let $g_{1},g_{2}:[0,a_{f})\rightarrow [0,\infty)$  be continuous functions such that for some $T\in (0,a_{f})$ we have $g_{2}(t)\leq g_{1}(t)$, $T\leq t<a_{f}$. If for a sequence $u_{m}$ of solutions of problem $N_{\lambda_{m}}$ we have
$$\int_{\Omega}g_{1}(u_{m})dx\leq C,$$
where $C$ is a constant independent of  $u_{m}$, then the same holds for $\int_{\Omega}g_{2}(u_{m})dx$.
\end{lem}
\begin{pf}
Indeed, we have
$$\int_{\Omega}g_{2}(u_{m})dx=\int_{u_{m}\leq T}g_{2}(u_{m})dx+\int_{u_{m}>T}g_{2}(u_{m})dx\leq M |\Omega|+\int_{\Omega}g_{1}(u_{m})dx$$
$$\leq M |\Omega|+C,~~\text{where}~~M:=\sup_{[0,T]}g_{2}(t).$$ $\blacksquare$
\end{pf}


\begin{prop}
Let $f$  satisfy $(H) $ (when $f$ is singular we additionally assume that $\lim_{t\rightarrow a_{f}}F(t)=\infty$). Let  $u_{m}$ be a sequence of semistable solutions of problem $(N_{\lambda_{m}})$. If  $\sup_{m}||\frac{\tilde{f}(u_{m})}{\sqrt{F(u_{m})}}||_{ L^{q}(\Omega)}<\infty$, for some $q\geq1$, then
\begin{equation}
\sup_{m}||u_{m}||_{L^{\infty}(\Omega)}<a_{f},
\end{equation}
for $n<2q$.
In particular, if  $\sup_{m}||f'(u_{m})||_{ L^{q}(\Omega)}<\infty$, then $(2.4)$ holds for $n<4q$.\\

\end{prop}
\begin{pf}
Take $v_{m}:=-\Delta u_{m}$, then from (1.1) $v_{m}$ satisfies
\begin{equation}
\left\{\begin{array}{ll} \Delta v_{m}+ \lambda_{m} f(u_{m})=0& {\rm }\ x\in \Omega,\\~~v_{m}=0& {\rm }\ x\in \partial \Omega.
\end{array}\right.
\end{equation}
We  rewrite problem (2.5) as $\Delta v_{m} + c(x)v_{m}=-\lambda_{m} f(0)$ where $c_{m}(x):=\lambda_{m} \frac{\tilde{f}(u_{m})}{v}$. By the pointwise estimate in Lemma 2.1 we have
$$0\leq c_{m}(x)=\lambda_{m} \frac{\tilde{f}(u_{m})}{v_{m}}\leq \sqrt{\lambda_{m}} \frac{\tilde{f}(u_{m})}{g(u_{m})},$$
Now using the inequality
$$0\leq\frac{\tilde{f}(t)}{g(t)}\leq \sqrt{2} \frac{\tilde{f}(t)}{\sqrt{F(t)}},  ~t>T,~~\text{for~some}~T<a_{f},$$
which comes from the fact that $\lim_{t\rightarrow a_{f}}\frac{g(t)}{\sqrt{2F(t)}}=1$,
and the assumptions, we get $\sup_{m}||c_{m}(x)||_{L^{q}(\Omega)}<\infty$. Thus, by the assumptions, Lemma 2.2 and Proposition 2.1, $||v_{m}||_{L^{\infty}(\Omega)}\leq C$, and hence  $||F( u_{m})||_{L^{\infty}(\Omega)}\leq C$  (by Lemma 2.1), where $C$ is a constant independent of $m$, for $n<2q$. Now the fact that $\lim_{t\rightarrow a_{f}}F(t)=\infty$ gives the first part. To prove the second part, it suffices to use Lemma 2.3 and  note that  by the convexity of $f$, we have
 \begin{equation}
\frac{\tilde{f}(t)}{\sqrt{F(t)}}\leq 2 \sqrt{f'(t)},~ \text{for}~ t ~\text{sufficiently close to}~ a_{f}.
\end{equation}
Indeed, $f'$ is a nondecreasing function by the convexity of $f$,  thus we have, for $0<t<a_{f}$
$$f'(t)F(t)=f'(t)\int_{0}^{t}f(s)ds\geq \int_{0}^{t}f'(s)f(s)ds=\frac{f(t)^{2}}{2}-\frac{f(0)^{2}}{2},$$
now the fact that $f(t)\rightarrow\infty$ as $t\rightarrow a_{f}$ gives (2.6).   $\blacksquare$
\end{pf}
 \begin{rem}
The condition $\lim_{t\rightarrow a_{f}}F(t)=\infty$ in the above proposition,   which is needed for a singular nonlinearity $f$, is satisfied by the extra assumption that $\tau_{+}<2$. Indeed, for a $\tau\in(\tau_{+},2)$ there exists $T\in(0,a_{f})$ such that $\frac{f''(t)}{f'(t)}\leq\tau \frac{f'(t)}{f(t)}$ for $t\in(T,a_{f})$, thus by an integration we get $f'(t)\leq Cf(t)^{\tau}$  or equivalently $f'(t)f(t)^{1-\tau}\leq Cf(t)$ for $t\in(T,a_{f})$. Again an integration gives
$$\frac{f(t)^{2-\tau}}{2-\tau}-\frac{f(T)^{2-\tau}}{2-\tau}\leq C(F(t)-F(T)),~~~\text{for}~~t\in(T,a_{f}).$$
Now the facts that $\lim_{t\rightarrow a_{f}}f(t)=\infty$ and $\tau<2$ imply that $\lim_{t\rightarrow a_{f}}F(t)=\infty$.
\end{rem}

For example, take the singular nonlinearity $f(t)=(1-t)^{-p}$ ($p>1$) on $[0,1)$. We have $\tau_{-}=\frac{p+1}{p}\in(0,2)$ and
$$F(t)=\frac{1}{p-1}\Big(\frac{1}{(1-t)^{p-1}}-1\Big)\rightarrow \infty,~~\text{as}~t\rightarrow1.$$
Then, as a corollary of Proposition 2.2 , we have the  next regularity result for problem ($N_{\lambda}$). It is proved in \cite{CG,CGE} by a different proof with the restriction that $p\neq 3$.
\begin{prop}
Let $f(u)=(1-u)^{-p}$ $(p>1)$ and $u_{m}$ be a sequence of semistable solutions of problem $(P_{\lambda_{m}})$, such that for some $q>1$ and $q\geq\frac{(p+1)n}{4p}$ so that $\sup_{m}||f(u_{m})||_{L^{q}(\Omega)}<\infty$. Then  $\sup_{m}||u_{m}||_{L^{\infty}(\Omega)}<1$.
\end{prop}
\begin{pf}
Notice that we have
$$f'(t)=p (1-t)^{-(p+1)}=f(t)^{\frac{p+1}{p}},~~t\in[0,1).$$
Hence, by the assumption $\sup_{m}||f'(u_{m})||_{L^{\frac{qp}{p+1}}(\Omega)}<\infty$.
\end{pf}

As an application of Proposition 2.2, consider problem $(N_{\lambda})$ with a convex nonlinearity $f$ satisfies ($H$) such that $f(t)=t\ln t$ for $t$ large. Then, for every $\epsilon>0$ there exist $T_{\epsilon},C_{\epsilon}>0$ such that $f'(t)\leq f(t)^{\epsilon}$ for $t\geq T_{\epsilon}$.
 \begin{equation}
f'(t)\leq f(t)^{\epsilon}~~\text{ for}~ t\geq T_{\epsilon}.
\end{equation}
Now if $u\geq0$ is a semistable solution of problem (1.1), from Lemma 3.5 \cite{CGE} we have $\int_{\Omega}f(u)dx\leq C$ with $C$ independent of $\lambda$ and $u$. This together (2.7) and Lemma 2.3 give $f'(u)\in L^{\frac{1}{\epsilon}}(\Omega)$ uniformly, hence by Proposition 2.2, $u^{*}$  is bounded for $n<\frac{4}{\epsilon}$, and since $\epsilon>0$ was arbitrary, $u^{*}$  is bounded in every dimension $n$. Indeed, the same result is true for every regular nonlinearity $f$ satisfies ($H$) with $\tau_{+}=0$ or equivalently
 \begin{equation}
\lim_{t\rightarrow\infty} \frac{f(t)f''(t)}{f'(t)^{2}}=0.
\end{equation}
Indeed, (2.8) implies (2.7) and we can proceed as above. \\

The following lemma is a special case of an interesting result of \cite{CGE}.
 \begin{lem}
Let $u$ be a semistable solution of problem $(N_{\lambda})$. If $H(t):=\int_{0}^{t}f''(s)\sqrt{F(s)}ds$ for $t\geq0$. Then
 \begin{equation}
\int_{\Omega}\sqrt{F(u)}H(u)dx\leq C,
\end{equation}
where $C$ is a constant independent of $\lambda$ and $u$.
\end{lem}

When $f$ is regular, in \cite{CGE} the authors used the above lemma to prove that  $u^{*}$  is bounded for $n<\frac{8}{\tau_{+}}$. In a completely similar manner and using Proposition 2.2, we can prove a similar result  when $f$ is singular.
\begin{lem}
Let $f$  satisfy $(H) $ be a singular nonlinearity with $0<\tau_{+}<2$, and  $u_{m}$ be a sequence of semistable solutions of problem $(N_{\lambda_{m}})$. Then
\begin{equation}
\sup_{m}||u_{m}||_{L^{\infty}(\Omega)}<a_{f},
\end{equation}
for $n<\frac{8}{\tau_{+}}$.
\end{lem}
\begin{pf}
 Take an arbitrary number $\tau>\tau_{+}$, then from the definition of $\tau_{+}$ there exists a $T_{1}\in(0,a_{f})$ such that  $\frac{f(t)f''(t)}{f'(t)^{2}}\leq \tau$, $T_{1}\leq t  <a_{f}$, which is equivalent to $\frac{d}{dt}(\frac{f'(t)}{\tilde{f}(t)^{\tau}})\leq0$ for $T_{1}\leq t  <a_{f}$. This gives $f'(t)\leq C_{0}f(t)^{\tau}$ for $T_{1}\leq t  <a_{f}$. Hence, using the inequality (2.6), $F(t)\geq C_{1}f'(t)^{\frac{2}{\tau}-1}$, $T_{1}\leq t  <a_{f}$, for some $T_{2}\in(T_{1},a_{f})$.
Thus, for a $T>T_{2}$ sufficiently close to $a_{f}$
we have
$$\sqrt{F(t)}H(t)\geq C_{2}f'(t)^{\frac{1}{\tau}-\frac{1}{2}} \int_{T}^{t}f''(s)f'(s)^{\frac{1}{\tau}-\frac{1}{2}} ds$$
$$\geq C_{3}f'(t)^{\frac{2}{\tau}},~ \text{for}~ t>T ~\text{sufficiently close to}~a_{f}.$$
Using the inequality above, Lemma 2.3 and Lemma 2.4,  we have $||f'(u_{m})||_{L^{\frac{2}{\tau}}}(\Omega) \leq C$. Hence by Remark 2.1  and Proposition 2.2, $\sup_{m}||u_{m}||_{L^{\infty}(\Omega)}<a_{f}$ for $n<\frac{8}{\tau}$, and since $\tau>\tau_{+}$ was arbitrary we get (2.10). $\blacksquare$
\end{pf}

As we have mentioned before, another main result of \cite{CGE} is that if $\tau_{-}>0$ then $u^{*}$ is bounded for $n\leq 7$. Using the same proof of this in \cite{CGE} we can prove it by a weaker assumption as follows:
\begin{cor}
Consider problem $(N_{\lambda})$ with a  regular nonlinearity $f$ satisfies  ($H$) such that for some $0\leq\epsilon<1$
 \begin{equation}
\liminf_{t\rightarrow\infty} \frac{f(t)^{1+\frac{\epsilon}{4}}f''(t)}{f'(t)^{2}}>0.
\end{equation}
Then $u^{*}$  is bounded for $n\leq7.$
\end{cor}
\begin{pf}
From (2.11) we have $f(t)^{1+\frac{\epsilon}{4}}f''(t)\geq C_{0}f'(t)^{2}$, $t\geq T$, for some $T>0$. Hence, using the inequality (2.6) and the fact that $F$ is a nondecreasing function we get, for a $T'>T$ sufficiently large,
$$\sqrt{F(t)}H(t)\geq C_{1}\int_{0}^{t}f''(s)F(s)ds\geq C_{2}\int_{T}^{t}\frac{f''(s)f(s)^{2}}{f'(s)}\geq C_{3}\int_{T}^{t}f'(s)f(s)^{1-\frac{\epsilon}{4}}$$
$$\geq C_{4}f(t)^{2-\frac{\epsilon}{4}}~ \text{for}~ t>T'.$$
Thus, from Lemmas 2.3 and 2.4 we have $||f(u)||_{L^{2-\frac{\epsilon}{4}}(\Omega)}<C$ where $C$ is independent of $u$. Now the elliptic regularity implies $u^{*}$  is bounded for $n\leq 8-\epsilon>7$, that gives the desired result. $\blacksquare$
\end{pf}


\section{Proof of the main results}

Following the idea of Dupaigne,  Ghergu and  Warnault in \cite{CG}, we prove the following lemma which is crucial for the proof of the main results.
\begin{lem}
Let $u$  be a positive smooth solution of $(N_{\lambda})$ satisfy the stability inequality $(1.2)$, $\theta:[0,a_{f})\rightarrow [0,\infty)$  a $C^{1}$ positive function with $\theta(0)=0$, and
$\Theta(t):=\int_{0}^{t}\theta'(s)^{2}ds$, for $0\leq t<a_{f}$. Then for every $\alpha>\frac{1}{2}$ we have
\begin{equation}
\int_{\Omega} \sqrt{f'(u)}\theta(u)^{2}dx\leq \frac{\alpha^{2}}{2\alpha-1} \Big(\int_{\Omega}     \frac{f(u)^{2\alpha}}{f'(u)^{\alpha-\frac{1}{2}}}                           dx\Big)^{\frac{1}{2\alpha}} \Big(\int_{\Omega} \frac{\Theta(u)^{\frac{2\alpha}{2\alpha-1}}}{f'(u)^{\frac{1}{2(2\alpha-1)}}} dx\Big)^{\frac{2\alpha-1}{2\alpha}}.
\end{equation}

\end{lem}
\begin{pf}
Let $u$  be a positive smooth solution of ($N_{\lambda}$) satisfy (1.2) and set $v:=-\Delta u$. Up to rescaling, we may assume that $\lambda=1$.  Take $\phi=\theta(u)$ as a test function in the stability inequality  (1.2). Then we get
\begin{equation}
\int_{\Omega} \sqrt{f'(u)}\theta(u)^{2}dx\leq \int_{\Omega} v\Theta(u)dx.
\end{equation}
Also, taking $\phi=v^{\alpha}$ ($\alpha>\frac{1}{2}$) as a test function in the stability inequality  (1.2), we get
\begin{equation}
\int_{\Omega} \sqrt{f'(u)}v^{2\alpha}dx\leq \frac{\alpha^{2}}{2\alpha-1}\int_{\Omega} f(u)v^{2\alpha-1}dx.
\end{equation}
Using H\"{o}lder inequality (with two conjugate numbers $2\alpha$ and $\frac{2\alpha}{2\alpha-1}$) on the right-hand side of inequality (3.2) we get
\begin{equation}
\int_{\Omega} \sqrt{f'(u)}\theta(u)^{2}dx\leq \Big(\int_{\Omega} \sqrt{f'(u)} v^{2\alpha}dx\Big)^{\frac{1}{2\alpha}} \Big(\int_{\Omega} \frac{\Theta(u)^{\frac{2\alpha}{2\alpha-1}}}{f'(u)^{\frac{1}{2(2\alpha-1)}}} dx\Big)^{\frac{2\alpha-1}{2\alpha}}.
\end{equation}
Similarly, form ( 3.3) and H\"{o}lder inequality we get
\begin{equation*}
\int_{\Omega} \sqrt{f'(u)}v^{2\alpha}dx\leq \frac{\alpha^{2}}{2\alpha-1}\Big(\int_{\Omega} \sqrt{f'(u)}v^{2\alpha}dx\Big)^{\frac{2\alpha-1}{2\alpha}}
\Big(\int_{\Omega}\frac{f(u)^{2\alpha}}{f'(u)^{\alpha-\frac{1}{2}}}dx\Big)^{\frac{1}{2\alpha}},
\end{equation*}
that gives
\begin{equation}
\int_{\Omega} \sqrt{f'(u)}v^{2\alpha}dx\leq \Big(\frac{\alpha^{2}}{2\alpha-1}\Big)^{2\alpha}
\int_{\Omega}     \frac{f(u)^{2\alpha}}{f'(u)^{\alpha-\frac{1}{2}}}                           dx.
\end{equation}
Plugging (3.5) in (3.4) we arrive at
\begin{equation*}
\int_{\Omega} \sqrt{f'(u)}\theta(u)^{2}dx\leq \frac{\alpha^{2}}{2\alpha-1} \Big(\int_{\Omega}     \frac{f(u)^{2\alpha}}{f'(u)^{\alpha-\frac{1}{2}}}                           dx\Big)^{\frac{1}{2\alpha}} \Big(\int_{\Omega} \frac{\Theta(u)^{\frac{2\alpha}{2\alpha-1}}}{f'(u)^{\frac{1}{2(2\alpha-1)}}} dx\Big)^{\frac{2\alpha-1}{2\alpha}},
\end{equation*}
which is the desired result. $\blacksquare$
\end{pf}
{\bf Proof of Theorem 1.1}\\
Fix an $\alpha>1$ such that $P_{f}(\alpha,\tau_{-},\tau_{+})<0$. Such an $\alpha$ exists since we have $P_{f}(1,\tau_{-},\tau_{+})=(2-\tau_{-})^{2}-4<0$ and $P_{f}(+\infty,\tau_{-},\tau_{+})=+\infty$ .
Now take positive numbers $\tau_{1}\in(0,\tau_{-})$  and $\tau_{2}\in(\tau_{+},2)$ such that
\begin{equation}
P_{f}(\alpha,\tau_{1},\tau_{2})<0.
\end{equation}
We claim that
\begin{equation}
I_{m}:=\int_{\Omega} \frac{\tilde{f}(u_{m})^{2\alpha}}{f'(u_{m})^{\alpha-\frac{1}{2}}}dx<C,
\end{equation}
where $C$ is independent of $m$. To this end, take $\theta(t)=\frac{\tilde{f}(t)^{\alpha}}{f'(t)^{\frac{\alpha}{2}}}$ in the inequality (3.1). First we estimate the function $\Theta(t)=\int_{0}^{t}\theta'(s)^{2}ds$ as follows. We have
\begin{equation}
\Theta(t)=\alpha^{2}\int_{0}^{t}\tilde{f}(s)^{2\alpha-2}f'(s)^{2-\alpha}\Big(1-\frac{\tilde{f}(s)f''(s)}{2f'(s)^{2}}\Big)^{2}ds.
\end{equation}
By the definitions of $\tau_{\pm}$ there exists a $T<a_{f}$ such that $\tau_{1}\leq\frac{\tilde{f}(t)f''(t)}{f'(t)^{2}}\leq \tau_{2}$ for $T\leq t <a_{f}$ that also gives
\begin{equation}
0<1-\frac{\tau_{2}}{2}\leq \frac{\tilde{f}(t)f''(t)}{2f'(t)^{2}}\leq 1-\frac{\tau_{1}}{2},~~\text{for}~ T\leq t <a_{f} .
\end{equation}
Using (3.9) in (3.8) we get
\begin{equation}
\Theta(t)\leq \Theta (T)+\alpha^{2} (1-\frac{\tau_{1}}{2})^{2}  \int_{T}^{t}\tilde{f}(s)^{2\alpha-2}f'(s)^{2-\alpha}ds,~~\text{for}~ T\leq t <a_{f} .
\end{equation}
Now, notice that taking $h(t):=\tilde{f}(t)^{2\alpha-1}f'(t)^{1-\alpha}$ for $0\leq t< a_{f}$, then
$$h'(t)=(2\alpha-1)\tilde{f}(t)^{2\alpha-2}f'(t)^{2-\alpha}\Big(1- \frac{\alpha-1}{2\alpha-1} \frac{\tilde{f}(s)f''(s)}{f'(s)^{2}}\Big)$$
$$\geq (2\alpha-1)(1-\frac{\alpha-1}{2\alpha-1}\tau_{2})\tilde{f}(t)^{2\alpha-2}f'(t)^{2-\alpha},~~\text{for}~ T\leq t <a_{f}.$$
Using the above inequality in (3.10) we obtain
\begin{equation}
\Theta(t)\leq C+ A  \tilde{f}(t)^{2\alpha-1}f'(t)^{1-\alpha},~\text{where}~A:=\frac{\alpha^{2}}{(2\alpha-1)} \frac{(1-\frac{\tau_{1}}{2})^{2}}{(1-\frac{\alpha-1}{2\alpha-1}\tau_{2})}~~\text{and}~C:=\Theta (T)-Ah(T).
\end{equation}
Note that in the above we also used that $1-\frac{\alpha-1}{2\alpha-1}\tau_{2}>0$ which holds since $\tau_{2}<2$.
Now, the fact that the inequality $\frac{\tilde{f}(t)f''(t)}{f'(t)^{2}}\leq \tau_{2}$ for $T\leq t <a_{f}$ is equivalent to $\frac{d}{dt}(\frac{f'(t)}{\tilde{f}(t)^{\tau_{2}}})\leq0$ for $T\leq t <a_{f}$  gives
\begin{equation}
f'(t)\leq C_{1}\tilde{f}(t)^{\tau_{2}}~~ \text{for}~ T\leq t <a_{f}.
\end{equation}
Using this we obtain, for $T\leq t <a_{f}$
$$\tilde{f}(t)^{2\alpha-1}f'(t)^{1-\alpha}\geq f'(t)^{\frac{2\alpha-1}{\tau_{2}}-(\alpha-1)}\rightarrow\infty,~\text{as}~t\rightarrow a_{f}.$$
Now take an $\epsilon>0$. From the inequality above and (3.11), there exists an $M_{\epsilon}\in[T,a_{f})$ such that
\begin{equation}
\Theta(t)\leq (A+\epsilon)\tilde{f}(t)^{2\alpha-1}f'(t)^{1-\alpha}, ~~\text{for}~ t\in[M_{\epsilon}, a_{f}).
\end{equation}
Hence,
\begin{equation}
\frac{\Theta(t)^{\frac{2\alpha}{2\alpha-1}}}{f'(t)^{\frac{1}{2(2\alpha-1)}}}\leq (A+\epsilon)^{\frac{2\alpha}{2\alpha-1}}\frac{\tilde{f}(t)^{2\alpha}}{f'(t)^{\alpha}}, ~~\text{for}~ t\in[M_{\epsilon}, a_{f}).
\end{equation}
Also, we can find an $M'_{\epsilon}>0$ such that
\begin{equation}
f(t)\leq(1+\epsilon)\tilde{f}(t),~~\text{for}~ t\in[M'_{\epsilon}, a_{f}).
\end{equation}
Now, taking $M''_{\epsilon}:=\max\{M_{\epsilon},M'_{\epsilon}\}$, then plugging (3.15), (3.14) in (3.1) we arrive at
$$I_{m}=
\int_{\Omega} \frac{\tilde{f}(u_{m})^{2\alpha}}{f'(u_{m})^{\alpha-\frac{1}{2}}}dx$$
$$\leq \frac{\alpha^{2}}{2\alpha-1} \Big(C_{\epsilon,m}+(1+\epsilon)^{2\alpha}\int_{u_{m}\geq M''_{\epsilon}}     \frac{\tilde{f}(u_{m})^{2\alpha}}{f'(u_{m})^{\alpha-\frac{1}{2}}}dx\Big)^{\frac{1}{2\alpha}} \Big(C'_{\epsilon,m}+(A+\epsilon)^{\frac{2\alpha}{2\alpha-1}}\int_{u_{m}\geq M''_{\epsilon}} \frac{\tilde{f}(u_{m})^{2\alpha}}{f'(u_{m})^{\alpha-\frac{1}{2}}} dx\Big)^{\frac{2\alpha-1}{2\alpha}},$$
where
$$C_{\epsilon,m}:=\int_{u<M''_{\epsilon}} \frac{f(u_{m})^{2\alpha}}{f'(u_{m})^{\alpha-\frac{1}{2}}}dx,~\text{and}~C'_{\epsilon,m}:=\int_{u_{m}<M''_{\epsilon}} \frac{\Theta(u_{m})^{\frac{2\alpha}{2\alpha-1}}}{f'(u_{m})^{\frac{1}{2(2\alpha-1)}}}dx.$$
Note that $C_{\epsilon,m}$ and $C'_{\epsilon,m}$ are bounded by a constant independent of $m$. Replacing the integrals on the right-hand side of the above inequality  with integrals over the full region $\Omega$ we get

\begin{equation}
I_{m}\leq \frac{\alpha^{2}}{2\alpha-1} \Big(C_{\epsilon,m}+(1+\epsilon)^{2\alpha}I_{m}\Big)^{\frac{1}{2\alpha}} \Big(C'_{\epsilon,m}+(A+\epsilon)^{\frac{2\alpha}{2\alpha-1}}I_{m}\Big)^{\frac{2\alpha-1}{2\alpha}}.
\end{equation}
Now if (3.7) does not hold then, $I_{m}\rightarrow\infty$ as $m\rightarrow\infty$. Hence, dividing two sides of  (3.16) by $I_{m}$ and letting $m\rightarrow\infty$, we must have
$$1\leq \frac{\alpha^{2}}{2\alpha-1}(1+\epsilon)(A+\epsilon),$$
and since $\epsilon>0$ was arbitrary we get
$$1\leq \frac{\alpha^{2}}{2\alpha-1}A=\frac{\alpha^{4}}{(2\alpha-1)^{2}}\frac{(1-\frac{\tau_{1}}{2})^{2}}{(1-\frac{\alpha-1}{2\alpha-1}\tau_{2})},$$
which is equivalent to $P_{f}(\alpha,\tau_{1},\tau_{2})\geq0$, a contradiction, that proves (3.7).\\
Now, using inequality (3.12) we have
$$\frac{\tilde{f}(t)^{2\alpha}}{f'(t)^{\alpha-\frac{1}{2}}}\geq C_{2}f'(t)^{\alpha(\frac{2}{\tau_{2}}-1)+\frac{1}{2}},~~for~ T\leq t<a_{f},$$
hence, thanks to Lemma 2.3 and (3.7) we get
$$||f'(u_{m})||_{L^{q_{1}}(\Omega)}\leq C,~~\text{where}~~q_{1}:=\alpha(\frac{2}{\tau_{2}}-1)+\frac{1}{2},$$
and $C$ is a constant independent of $m$. Now Proposition 2.2 implies that
\begin{equation}
\sup_{m}||u_{m}||_{L^{\infty}(\Omega)}<a_{f},~~for~~n<4\alpha(\frac{2}{\tau_{2}}-1)+2.
\end{equation}
Again from inequality (3.12) we have
$$\frac{\tilde{f}(t)^{2\alpha}}{f'(t)^{\alpha-\frac{1}{2}}}\geq C_{3}f(t)^{\alpha(2-\tau_{2})+\frac{\tau_{2}}{2}},~~for~T\leq t<a_{f},$$
and using the above  inequality, Lemma 2.3 and  (3.7) we get
$$||f(u_{m})||_{L^{q_{2}}(\Omega)}\leq C,~~\text{where}~~q_{2}:=\alpha(2-\tau_{2})+\frac{\tau_{2}}{2},$$
and $C$ is a constant independent of $m$.
Hence, in the case when $f$ is regular, by the elliptic regularity theory
\begin{equation}
\sup_{m}||u_{m}||_{L^{\infty}(\Omega)}<\infty,~~for~~n<4\alpha(2-\tau_{2})+2\tau_{2}.
\end{equation}
Now, since we can choose $\tau_{2}$ arbitrary close to $\tau_{+}$ and $\alpha$ near to the largest root of the polynomial $P_{f}$, then  (3.17) and  (3.18) complete the proof of the first part.\\
To see the second part, suppose that $\tau_{-}=\tau_{+}:=\tau>0$. If $\tau<\frac{2}{3}$ then  from Lemma 2.5, $\sup_{m} ||u_{m}||_{L^{\infty}(\Omega)}<a_{f}$ for $n\leq \frac{8}{\tau}>12$, so we need to prove it for the case $\frac{2}{3}\leq\tau \leq1$. It is not hard to see (for example by using a computing device) that, for $\alpha=\frac{5\tau}{2(2-\tau)}$ we have $P_{f}(\alpha,\tau,\tau)<0$ on the interval $[\frac{2}{3},1]$, hence $\alpha^{*}>\frac{5\tau}{2(2-\tau)}$ that gives
$$\sup_{m} ||u_{m}||_{L^{\infty}(\Omega)}<a_{f},~~for~~n<4\alpha^{*}(\frac{2}{\tau}-1)+2>12.$$
Also, when $1\leq \tau\leq1.57863$ then for $\alpha=\frac{5\tau}{4(2-\tau)}$ we have $P_{f}(\alpha,\tau,\tau)<0$ on the interval $[\frac{2}{3},1]$, hence $\alpha^{*}>\frac{5\tau}{4(2-\tau)}$ that gives
$$\sup_{m} ||u_{m}||_{L^{\infty}(\Omega)}<a_{f},~~for~~n<4\alpha^{*}(\frac{2}{\tau}-1)+2>7,$$
and now the proof is complete. $\blacksquare$\\

\section{Acknowledgement}
This research was in part supported by a grant from IPM (No. 94340123).



\begin{thebibliography}{}

\bibitem{ADN} S. Agmon, A. Douglis, L. Nirenberg, Estimates near the boundary for solutions of elliptic partial differential equations satisfying general boundary conditions, I, Comm. Pure Appl. Math. 12 (1959) 623–727.
 \bibitem{BG} E. Berchio, F. Gazoola, Some remarks on bihormonic elliptic problems with positive, increasing and convex nonlinearities, Electronic J. differential Equations, 34 (2005) 20 pp.


\bibitem{CaG} D. Cassani, J. do O, N. Ghoussoub, On a fourth order elliptic problem with a singular nonlinearity, Adv. Nonlinear Stud. 9 (2009) 177-197.
 \bibitem{C}    C. Cowan, Regularity of the extremal solutions in a Gelfand system problem, Adv. Nonlinear
Stud., 11 (2011), 695–700.

 \bibitem{CGEM} C. Cowan, P. Esposito, N. Ghoussoub and A. Moradifam, The critical dimension for a
fourth order  elliptic problem with singular nonlinearity, Arch. Ration. Mech. Anal., in
press (2009) 19 pp.
\bibitem{CGE} C. Cowan, P. Esposito and N. Ghoussoub, Regularity of extremal solutions in fourth order nonlinear eigenvalue problems on general domains, DCDS-A 28 (2010), 1033-1050.
\bibitem{CG} C. Cowan, N. Ghoussoub, Regularity of semi-stable solutions to fourth order nonlinear eigenvalue problems on general domains, Cal. Var. 49 (2014) 291-305.
\bibitem{CSS} X. Cabr\'{e}, M. Sanch\'{o}n and J. Spruck, A priori estimates for semistable solutions of semilinear elliptic equations, arXiv preprint, arXiv:1407.0243, 2014 - arxiv.org (2014).
 \bibitem{DDG}    J. D´avila, L. Dupaigne, I. Guerra and M. Montenegro, Stable solutions
for the bilaplacian with exponential nonlinearity, SIAM J. Math. Anal.
39 (2007), 565-592.


  \bibitem{DFG}  J. D´avila, I. Flores and I. Guerra, Multiplicity of solutions for a fourth
order equation with power-type nonlinearity, Math. Ann. 348 (2010),
143-193.





 \bibitem{Dup1} L. Dupaigne, M. Ghergu and G. Warnault,  The Gel'fand Problem for the Biharmonic Operator, Arch. Ration. Mech. Anal., 208, 3 (2013) 725-752.

    \bibitem{Dup2}  L. Dupaigne, A. Farina, and B. Sirakov, Regularity of the extremal solution
for the Liouville system, pp. 139-144 in Geometric partial differential equations (Pisa,
2012), edited by A. Chambolle et al., CRM Series 15, Ed. Norm., Pisa, Pisa, 2013. MR 3156892
arXiv 1207.3703v1.


     \bibitem{Fer} A. Ferrero, H.-C. Grunau, and P. Karageorgis, Supercritical biharmonic
equations with power-type nonlinearity, Ann. Mat. Pura Appl. (4) 188:1 (2009) 171–185.

\bibitem{GuW} Z Guo, J Wei, Liouville type results and regularity of the extremal solutions of biharmonic equation with negative exponents, Discrete Contin. Dyn. Syst,  34 (2014), 2561-2580.

   \bibitem{Amir}  A. Moradifam, The singular extremal solutions of the bilaplacian with exponential nonlinearity, Proc. Amer. Math. Soc. 138 (2010), 1287-1293.
       \bibitem{Gaz} F. Gazzola, H.-C. Grunau, and G. Sweers, Polyharmonic boundary value
problems: positivity preserving and nonlinear higher order elliptic equations in bounded domains,
Lecture Notes in Mathematics 1991, Springer, Berlin, 2010.

 \bibitem{GW}   Z . Guo. J. Wei, On a fourth order nonlinear elliptic equation with negative exponent,
SIAM J. Math. Anal. 40 (2008/09), 2034-2054.


\bibitem{HHY} H. Hajlaoui, A. Harrabi, D. Ye, On stable solutions of the biharmonic problem with polynomial growth, Pacific Journal of Mathematics 270, 1 (2014) 79-93.
    \bibitem{Ser} J. Serrin, Local behavior of solutions of quasi-linear equations, Acta Math. 111
(1964) 247-302.
    \bibitem{W} K. Wang, Partial regularity of stable solutions to the supercritical equations and its applications,
Nonlinear Anal., 75 (2012), 5238–5260.



    \bibitem{WXY} J. Wei, X. Xu, and W. Yang, "On the classification of stable solutions to biharmonic
problems in large dimensions", Pacific J. Math. 263:2 (2013), 495–512. MR 3068555 Zbl 06196725

\bibitem{YW} D. Ye, J. Wei, Liouville Theorems for finite Morse index solutions of
Biharmonic problem, Math. Ann. to appear.



\end{thebibliography}
\end{document}